\documentclass[12pt,a4paper,reqno]{amsart}
\usepackage{amssymb}
\usepackage{dcpic,pictex}
\usepackage{multirow}
\usepackage{graphicx}
\usepackage{slashbox}
\usepackage{cite}

\usepackage{hyperref}

\newtheorem{The}{Theorem}[section]
\newtheorem{Rem}{Remark}[section]

\newcommand{\R}{\mathbb{R}}
\newcommand{\C}{\mathbb{C}}
\newcommand{\N}{\mathbb{N}}

\begin{document}


\title[On complete monotonicity]{On complete monotonicity of solution to the fractional relaxation equation with the $n$th level fractional derivative}

\author{Yuri Luchko}
\curraddr{Beuth Technical University of Applied Sciences Berlin,  
     Department of  Mathematics, Physics, and Chemistry,  
     Luxemburger Str. 10,  
     13353 Berlin,
Germany}
\email{luchko@beuth-hochschule.de}


\begin{abstract}
In this paper, we first deduce the explicit formulas for the projector of the $n$th level fractional derivative and for its Laplace transform. Then the fractional relaxation equation with the $n$th level fractional derivative is discussed. It turns out that under some conditions, the solutions to the initial-value problems for this equation are completely monotone functions that can be represented in form of the linear combinations of the Mittag-Leffler functions with some power law weights. Special attention is given to the case of the relaxation equation with the 2nd level derivative.
\end{abstract}

\maketitle

\section{Introduction}
\label{s1}

In the framework of the current discussions regarding the ``right" fractional derivatives (\cite{DGGS, Giu18, Han20, Hil19, HilLuc, KocLuc19, Luc20, Ort15, Ort18, Sty18, Tar13, Tar19}), the main suggested approach was to define the fractional derivatives via the Fundamental Theorem of Fractional Calculus (FC), i.e., as the left-inverse operators to the corresponding fractional integrals that satisfy the index low, interpolate the definite integral, and build a family of the operators continuous in a certain sense with respect to the order of integration. 

According to a result derived in \cite{Cart78}, under the conditions mentioned above, the only family of the fractional integrals defined on a finite interval are the Riemann-Liouville fractional integrals \cite{Samko}. Until recently, three families of the fractional derivatives that are the left-inverse operators to the Riemann-Liouville fractional integrals were discussed in the literature: the Riemann-Liouville fractional derivatives \cite{Samko}, the Caputo fractional derivatives \cite{Diet10}, and the Hilfer fractional derivatives \cite{Hil00}.  However, in \cite{Luc20}, infinitely many other families of the fractional derivatives that are the left-inverse operators to the Riemann-Liouville fractional integrals, were introduced and called the $n$th level fractional derivatives. These derivatives satisfy the Fundamental Theorem of FC, i.e., they are the left-inverse operators to the Riemann-Liouville fractional integrals on the appropriate nontrivial spaces of functions that justifies calling them the fractional derivatives. 

In \cite{Luc20}, some basic properties of the $n$th level fractional derivatives were studied including a description of their kernels. However, the question of their applicability to some real world problems remained open. In this paper, we provide a first evidence of their usefulness for applications on an example from the linear viscoelasticity. More precisely, we show that the  solution to the Cauchy problem for the fractional relaxation equation with the $n$th level fractional derivative is a completely monotone function that can be represented in form of a linear combination of the Mittag-Leffler functions with some power law weights. 
As discussed in \cite{Mai10}, the property of complete monotonicity  is characteristic for any relaxation process. Only in this case it can be interpreted as a superposition of (infinitely many) elementary, i.e., exponential, relaxation processes. In linear viscoelasticity, the assumption of complete monotonicity of the solutions to the relaxation equations that model the relaxation processes is usually supposed to be fulfilled, see, e.g., \cite{Mai10} and references therein.

The rest of this paper is organized as follows: in Section 2, we discuss some new properties of the  $n$th level fractional derivative including the explicit formulas for the projector  of the $n$th level fractional derivative and for its Laplace transform.  Section 3 addresses the Cauchy problem for the relaxation equation with the  $n$th level fractional derivative and properties of its solution including complete monotonicity. In particular, we focus on an important particular case of the Cauchy problem for the relaxation equation with the  2nd level fractional derivative.

\section{The $n$th level fractional derivative and its properties}

\subsection{Basic definitions and properties}

In \cite{Luc20}, the $n$th level fractional derivative on a finite interval (without loss of generality we proceed with the interval $[0,\, 1]$) was introduced  as follows:
let  $0<\alpha \le 1$ and the parameters $\gamma_1,\ \gamma_2,\dots,\gamma_n \in \R$ satisfy the conditions 
\begin{equation}
\label{gamma}
0\le \gamma_k \ \mbox{and}\  \alpha + s_k \le k, \ s_k:= \sum_{i=1}^k\, \gamma_i,\ k=1,2,\dots,n.
\end{equation}
The operator 
\begin{equation}
\label{nLD}
(D^{\alpha,(\gamma)}_{nL}\, f)(x)  = \left(\prod_{k=1}^n (I^{\gamma_k}\, \frac{d}{dx})\right)\, (I^{n-\alpha-s_n}\, f)(x)
\end{equation}
is called the $n$th level  fractional derivative of order  $\alpha$ and type $\gamma = (\gamma_1,\ \gamma_2,\dots,\gamma_n)$.

In \eqref{nLD}, $I^\alpha$ stands for the Riemann-Liouville fractional integral ($x\in[0,\, 1]$):
\begin{equation}
\label{RLI}
(I^\alpha\, f)(x) = \begin{cases}
\frac{1}{\Gamma(\alpha)} \int_0^x (x-t)^{\alpha -1}\, f(t)\, dt, & \alpha >0,\\
f(x), & \alpha = 0
\end{cases}
\end{equation}

The $n$th level  fractional derivative  is well defined, say, on the function space
\begin{equation}
\label{X1nL}
X_{nL} = \{ f: \left(\prod_{k=i}^n (I^{\gamma_k}\, \frac{d}{dx})\right)\, I^{n-\alpha-s_n}\, f \in \mbox{AC}([0,\, 1]),\ i=2,\dots n+1 \}, 
\end{equation} 
where the notation $\mbox{AC}([0,1])$ stands for the space of the absolutely continuous functions (an empty product is interpreted as the identity operator). The absolutely continuous functions allow the following representation:
\begin{equation}
\label{absc}
f\in \mbox{AC}([0,1]) \, \Leftrightarrow \, \exists \phi \in L_1(0,1):\, f(x) = f(0) + \int_0^x \phi(t)\, dt,\ x\in [0,\, 1].
\end{equation}
In what follows, a (weak) derivative of a function $f\in \mbox{AC}([0,1])$ will be understood in the following sense:
\begin{equation}
\label{absc_d}
f(x) = f(0) + \int_0^x \phi(t)\, dt,\ x\in [0,\, 1] \ \Rightarrow \ \frac{df}{dx} := \phi \in L_1(0,1).
\end{equation}
Then the $n$th level  fractional derivative is a linear operator that maps $X_{nL}$ into $L_1(0,1)$ (see \cite{Luc20}).

The notation "$n$th level  fractional derivative" is justified by the Fundamental Theorem of FC.

\begin{The}[\cite{Luc20}]
\label{t_ft_Ln}
Let $X_{FT}$ be the following space of functions:
\begin{equation}
\label{X2}
X_{FT} = \left\{ f:\, I^{\alpha}f\in \mbox{AC}([0,\, 1])\ \mbox{and}\ (I^{\alpha}f)(0) = 0\right\}.
\end{equation}
The $n$th level  fractional derivative is a left-inverse operator to the Riemann-Liouville fractional integral on the space $X_{FT}$  , i.e., the relation
\begin{equation}
\label{ftFC_nL}
(D^{\alpha,(\gamma)}_{nL}\, I^\alpha\, f)(x) = f(x),\ x\in [0,\, 1]
\end{equation} 
holds true for any $f$ from the space $X_{FT}$.
\end{The}

\begin{Rem}
\label{r_ft}
In calculus, the formula of type \eqref{ftFC_nL} $\left(\frac{d}{dx}\, \int_0^x\, f(t)\, dt\ =\, f(x)\right)$ is usually called the 1st fundamental theorem of calculus. The 2nd fundamental theorem of calculus states that $\int_{0}^x\, f^\prime(t)\, dt \, = \, f(x)-f(0)$. We address the 2nd fundamental theorem of FC for the $n$th level  fractional derivative in the next subsection. 
\end{Rem}

\begin{Rem}
\label{r_dn}
In \cite{DN}, the differential operators similar to the $n$th level fractional derivative $D^{\alpha,(\gamma)}_{nL}$  (in other notations and with other restrictions on the parameters) have been addressed. 

However, the motivation behind the $n$th level fractional derivative is essentially different compared to the one behind the Djrbashian-Nersesian differential operator introduced in \cite{DN}. The main feature of the $n$th level fractional derivative is its connection to the Riemann-Liouville fractional integral of order $\alpha$ (Theorem \ref{t_ft_Ln}).  The parameter $\alpha$ in its definition plays the role of the order of this fractional derivative. Moreover, the restriction $0<\alpha \le 1$ is very important in this context, because only in this case the $n$th level fractional derivative in its present form is a left inverse operator to the Riemann-Liouville fractional integral (see the proof of Theorem \ref{t_ft_Ln} in \cite{Luc20}). 

In \cite{DN}, this connection was not discussed at all. The type $(\gamma_0,\gamma_1,\dots,\gamma_n)$ of the operators from \cite{DN} does not have an explicit connection to the Riemann-Liouville integral because all parameters are ``uniform", there is no selected parameter that stands for the order of this operator. Another difference is the restrictions on the parameters for the $n$th level fractional derivative and for the operators introduced in \cite{DN}. They are essentially different because the $n$th level fractional derivative is a fractional derivative in the sense of the Fundamental Theorem of FC, whereas the operators from \cite{DN} are just some general ``differential operations". Finally, the Djrbashian-Nersesian operator was introduced and employed only in the case of a finite interval, whereas in this paper, the $n$th level fractional derivative is addressed also on the positive real semi-axis.
\end{Rem}

Because of the fundamental theorem of calculus  and the index law for the Riemann-Liouville fractional integrals $\left( I^{\alpha +\beta} \, = \, I^\alpha\, I^\beta,\ \alpha,\beta \ge 0\right)$, the $n$th level fractional derivatives are reduced to the fractional derivatives of the level less than $n$ if  some of the parameters $\gamma_k,\ k=2,\dots,n$ are equal to or grater than one or if the inequality $\alpha+s_n\le n-1$ holds valid. That's why, in the following discussions, we address only the truly $n$th level fractional derivatives and suppose that the conditions
\begin{equation}
\label{cond_add_n}
n-1 <\alpha+s_n \ \mbox{and}\  \gamma_k <1,\ k= 2,\dots,n
\end{equation} 
hold valid.

Under these conditions, the kernel of the $n$th level fractional derivative is $n$-dimensional (\cite{Luc20}):
\begin{equation}
\label{nLD_K}
\mbox{Ker}(D^{\alpha,(\gamma)}_{nL}) =\left\{\sum_{k=1}^n c_k x^{\sigma_k},\ \sigma_k= \alpha +s_k -k,\ c_k\in \R \right\}.
\end{equation}
Because of the conditions \eqref{gamma} and \eqref{cond_add_n}, the exponents $\sigma_k$ of the basis functions of the kernel fulfill the inequalities 
\begin{equation}
\label{sigma}
-1 <\sigma_k\le 0,\ k=1,2,\dots,n.
\end{equation}
In the case, one or several of the conditions from \eqref{cond_add_n} do not hold true, the $n$th level fractional derivatives degenerate to the derivatives of the level less than $n$ and thus their kernels have dimensions less 
than $n$. 

\subsection{Projector of the $n$th level fractional derivative}

One of the most important and widely used methods for analysis of the fractional differential equations is by means of their reduction to certain integral equations of Volterra type. To perform this transformation, one acts on the fractional differential equations with the corresponding fractional integrals. Thus, for this method, the explicit formulas for the compositions of the fractional integrals and the fractional derivatives (2nd Fundamental Theorem of FC) are required. In this subsection, we derive a formula of this kind for the  $n$th level fractional derivative. 

\begin{The}
\label{t-pro}
Under the conditions \eqref{cond_add_n},  the projector  
\begin{equation}
\label{pro}
(P^\alpha_{nL}\, f)(x) = (\mbox{Id}-I^{\alpha} D^{\alpha,(\gamma)}_{nL}\, f)(x)
\end{equation}
of the $n$th level fractional derivative \eqref{nLD}  on the function space $X_{nL}$ defined by \eqref{X1nL} takes the form
\begin{equation}
\label{pro1}
(P^\alpha_{nL}\, f)(x) = \sum_{k=1}^n \, p_k\, x^{\sigma_k}, \ \sigma_k = \alpha +s_k -k,
\end{equation}
\begin{equation}
\label{p}
p_k = \frac{\left( \prod_{i=k+1}^n \left(I^{\gamma_i}\, \frac{d}{dx}\right) \, I^{n-\alpha-s_n}\, f\right)(0)}{\Gamma(\alpha +s_k -k +1)}.
\end{equation}
\end{The}

We start the proof of the theorem by introducing an auxiliary function
\begin{equation}
\label{g}
g(x) := (I^{\alpha} D^{\alpha,(\gamma)}_{nL}\, f)(x).
\end{equation}
For $f\in X_{nL}$, the derivative $D^{\alpha,(\gamma)}_{nL}\, f$ is well defined and belongs to $L_1(0,1)$. Thus the function $g$ is from the space $I^\alpha(L_1(0,1))$. Then we can act on $g$ with  the operator $D^{\alpha,(\gamma)}_{nL}$ and apply the Fundamental Theorem of FC for the $n$th level fractional derivative:
$$
(D^{\alpha,(\gamma)}_{nL}\, g)(x) = (D^{\alpha,(\gamma)}_{nL}\, I^{\alpha} D^{\alpha,(\gamma)}_{nL}\, f)(x) = (D^{\alpha,(\gamma)}_{nL}\, f)(x).
$$
This formula means that the function $g-f$ belongs to the kernel of the $n$th level fractional derivative given by the formula \eqref{nLD_K}. Thus, we get the representation 
\begin{equation}
\label{f-g}
g(x) = f(x) + \sum_{k=1}^n \, c_k\, x^{\sigma_k},\ \sigma_k=\alpha +s_k -k.
\end{equation} 

Now we determine the coefficients $c_k,\ k=1,\dots,n$. To do this, we apply a sequence of some Riemann-Liouville fractional integrals to the auxiliary function $g$ defined by \eqref{g} and employ the 1st fundamental theorem of calculus and the well-known formula
\begin{equation}
\label{power}
(I^\alpha\, t^{\beta})(x) = \frac{\Gamma(\beta +1)}{\Gamma(\alpha+\beta +1)}\, x^{\alpha + \beta},\ \alpha \ge 0,\ \beta > -1.
\end{equation} 

First we apply the Riemann-Liouville fractional integral $I^{1 -\alpha-\gamma_1}$:
$$
(I^{1-\alpha-\gamma_1 }\, g)(x) = (I^{1-\alpha-\gamma_1 } I^{\alpha} D^{\alpha,(\gamma)}_{nL}\, f)(x) =
$$
$$
\left(I^{1-\alpha-\gamma_1 } I^{\alpha} \left(\prod_{k=1}^n (I^{\gamma_k}\, \frac{d}{dx})\right)\, (I^{n-\alpha-s_n}\, f)\right)(x) =
$$
$$
\left(I^{1-\gamma_1 }  I^{\gamma_1} \frac{d}{dx} \left(\prod_{k=2}^n (I^{\gamma_k}\, \frac{d}{dx})\right)\, (I^{n-\alpha-s_n}\, f)\right)(x) =
$$
$$
\left(I^{1}\, \frac{d}{dx}\, \left(\prod_{k=2}^n (I^{\gamma_k}\, \frac{d}{dx})\right)\, (I^{n-\alpha-s_n}\, f)\right)(x) =
$$
$$
\left(\left(\prod_{k=2}^n (I^{\gamma_k}\, \frac{d}{dx})\right)\, (I^{n-\alpha-s_n}\, f)\right)(x) - \left(\left(\prod_{k=2}^n (I^{\gamma_k}\, \frac{d}{dx})\right)\, (I^{n-\alpha-s_n}\, f)\right)(0).
$$
Then we apply the Riemann-Liouville fractional integral  $I^{1-\gamma_2}$: 
$$
(I^{1-\gamma_2}\, I^{1 -\alpha-\gamma_1}\, g)(x) = \left(I^{1}\, \frac{d}{dx}\, \left(\left(\prod_{k=3}^n (I^{\gamma_k}\, \frac{d}{dx})\right)\, (I^{n-\alpha-s_n}\, f)\right)\right)(x) - 
$$
$$
\frac{\left(\left(\prod_{k=2}^n (I^{\gamma_k}\, \frac{d}{dx})\right)\, (I^{n-\alpha-s_n}\, f)\right)(0)}{\Gamma(2-\gamma_2)}\, x^{1-\gamma_2} =
\left(\left(\prod_{k=3}^n (I^{\gamma_k}\, \frac{d}{dx})\right)\, (I^{n-\alpha-s_n}\, f)\right)(x) - 
$$
$$
\left(\left(\prod_{k=3}^n (I^{\gamma_k}\, \frac{d}{dx})\right)\, (I^{n-\alpha-s_n}\, f)\right)(0)- 
\frac{\left(\left(\prod_{k=2}^n (I^{\gamma_k}\, \frac{d}{dx})\right)\, (I^{n-\alpha-s_n}\, f)\right)(0)}{\Gamma(2-\gamma_2)}\, x^{1-\gamma_2}.
$$
Now we continue by applying the Riemann-Liouville fractional integrals $I^{1-\gamma_k},\ k=3,\dots,n$ and thus arrive at the following formula:
$$
\left(\left(\prod_{k=2}^n I^{1-\gamma_k}\right) (I^{1 -\alpha-\gamma_1}\, g)\right)(x) = (I^{n-\alpha - s_n} g)(x) = 
$$
$$
(I^{n-\alpha - s_n} f)(x) - \sum_{k=1}^n \frac{\left(\left(\prod_{i=k+1}^n (I^{\gamma_k}\, \frac{d}{dx})\right)\, (I^{n-\alpha-s_n}\, f)\right)(0)}{\Gamma(n-k+1-(s_n - s_k))}\, x^{n-k-(s_n - s_k)},
$$
where the empty product is understood as the identity operator. 
On the other hand, we can apply the operator $\left(\prod_{k=2}^n I^{1-\gamma_k}\right) I^{1 -\alpha-\gamma_1} = I^{n-\alpha - s_n} $ to the identity \eqref{f-g} and thus get the relation 
$$
(I^{n-\alpha - s_n} g)(x) = (I^{n-\alpha - s_n} f)(x) + 
\sum_{k=1}^n c_k \frac{\Gamma(\alpha+s_k - k+1))}{\Gamma(n-k+1-(s_n - s_k))}\, x^{n-k-(s_n - s_k)}.
$$
The coefficients $c_k,\ k=1,\dots,n$ are obtained by comparison of the coefficients by the same powers of $x$ in the last two formulas:
$$
c_k = -\frac{\left( \prod_{i=k+1}^n \left(I^{\gamma_i}\, \frac{d}{dx}\right) \, I^{n-\alpha-s_n}\, f\right)(0)}{\Gamma(\alpha +s_k -k +1)},\ k=1,\dots,n.
$$
The statement of the theorem follows from this formula and the representation \eqref{f-g}.

\begin{Rem}
\label{r_2ndFT}
Theorem \ref{t-pro} can be rewritten in form of the 2nd Fundamental Theorem of FC for the $n$th level fractional derivative:
\begin{equation}
\label{2ndFT}
(I^{\alpha} D^{\alpha,(\gamma)}_{nL}\, f)(x) = f(x) - \sum_{k=1}^n \, p_k\, x^{\sigma_k}, \ \sigma_k = \alpha +s_k -k,
\end{equation}
where the coefficients $p_k,\ k=1,\dots,n$ depend on the function $f$ and are given by the formula  \eqref{p}.

If some of the conditions \eqref{cond_add_n} do not hold true, the $n$th level fractional derivatives are reduced to the fractional derivatives with the level less than $n$ and their kernels have the dimensions less than $n$. In this case, some of the coefficients $c_k$ in the representation \eqref{f-g} and thus the corresponding coefficients $p_k$ in the formula \eqref{pro1} for the projector 
$P^\alpha_{nL}$ are equal to zero.
\end{Rem}

As an example, let us discuss the projector of the 2nd level fractional derivative deduced in \cite{Luc20}. In the general case, it has the following form:
\begin{equation}
\label{pro1_2}
(P^\alpha_{2L}\, f)(x) = p_1\, x^{\alpha +\gamma_1-1} + p_2\, x^{\alpha +\gamma_1+\gamma_2-2},
\end{equation}
\begin{equation}
\label{p1_2}
\begin{cases}
p_1 = \frac{1}{\Gamma(\alpha +\gamma_1)}\left( I^{\gamma_2}\, \frac{d}{dx} \, I^{2-\alpha-\gamma_1-\gamma_2}\, f\right)(0), \\ 
p_2= \frac{1}{\Gamma(\alpha +\gamma_1+\gamma_2 -1)}\left(  I^{2-\alpha-\gamma_1-\gamma_2}\, f\right)(0).
\end{cases}
\end{equation}
If one of the conditions 
\begin{equation}
\label{cond_add}
\gamma_2 <1,\ 1<\alpha+\gamma_1+\gamma_2
\end{equation} 
does not hold true, the 2nd level fractional derivative is reduced to the Hilfer fractional derivative defined by
\begin{equation}
\label{HD}
(D^{\alpha,\gamma_1}_H\, f)(x)  = (I^{\gamma_1}\, \frac{d}{dx}\, I^{1-\alpha-\gamma_1}\, f)(x).
\end{equation}
Indeed, for $1 \le \gamma_2$, we get the following chain of equalities:
$$
D^{\alpha,(\gamma_1,\gamma_2)}_{2L} = I^{\gamma_1}\, \frac{d}{dx}\, I^1\,I^{\gamma_2-1} \, \frac{d}{dx}\,I^{2-\alpha-\gamma_1-\gamma_2} = 
I^{\gamma_1+\gamma_2 -1}\, \frac{d}{dx}\, I^{1-\alpha-(\gamma_1+\gamma_2-1)} = D^{\alpha,\gamma_1+\gamma_2-1}_{H}. 
$$
In the case $\alpha +\gamma_1+\gamma_2 \le 1 $, we proceed as before:  
$$
D^{\alpha,(\gamma_1,\gamma_2)}_{2L} = I^{\gamma_1}\, \frac{d}{dx}\, I^{\gamma_2} \, \frac{d}{dx}\,I^{1}\, I^{1-\alpha-\gamma_1-\gamma_2}
=  
I^{\gamma_1}\, \frac{d}{dx}\, I^{\gamma_2} I^{1-\alpha-\gamma_1-\gamma_2} = 
 I^{\gamma_1}\, \frac{d}{dx}\,  I^{1-\alpha-\gamma_1} = D^{\alpha,\gamma_1}_{H}.
$$ 
In both cases, the kernels of the 2nd level fractional derivatives (the Hilfer fractional derivatives) are one-dimensional and one of the coefficients $c_1$ or $c_2$ in the representation \eqref{f-g} and thus one of the coefficients $p_1$ or $p_2$ in the formula \eqref{pro1} for the projector  
$P^\alpha_{2L}$ is equal to zero. As a result, the projector $P_H^\alpha$ of the Hilfer fractional derivative  takes the following known form (\cite{Hil}):
\begin{equation}
\label{pro1_H_1}
(P^\alpha_{H}\, f)(x) = \frac{1}{\Gamma(\alpha +\gamma_1)}\left(  I^{1-\alpha-\gamma_1}\, f\right)(0)\, x^{\alpha +\gamma_1-1}.
\end{equation}
Substituting $\gamma_1 = 0$ into the formula \eqref{pro1_H_1}, we get the projector of the Riemann-Liouville fractional derivative
\begin{equation}
\label{pro1_RL}
(P^\alpha_{RL}\, f)(x) = \frac{1}{\Gamma(\alpha)}\left(  I^{1-\alpha}\, f\right)(0)\, x^{\alpha -1}.
\end{equation} 
The value $\gamma_1 = 1-\alpha$ corresponds to the projector of the Caputo fractional derivative:
\begin{equation}
\label{pro1_C}
(P^\alpha_{C}\, f)(x) = f(0).
\end{equation}
The Riemann-Liouville and the Caputo fractional derivatives are both particular cases of the Hilfer fractional derivative and are defined as follows:
\begin{equation}
\label{RLD}
(D^\alpha_{RL}\, f)(x) = \frac{d}{dx}(I^{1-\alpha} f)(x),
\end{equation}
\begin{equation}
\label{CD}
(D^\alpha_{C}\, f)(x) = (I^{1-\alpha} \frac{d}{dx}f)(x).
\end{equation}

\subsection{Laplace transform of the $n$th level fractional derivative}

In this subsection, we introduce the $n$th level fractional derivative on the positive real semi-axis and derive a formula for its Laplace transform. 

The Riemann-Liouville fractional integral on the positive real semi-axis is given by the formula \eqref{RLI} extended to the case $x\in \R_+$:
\begin{equation}
\label{RLI_R}
(I^\alpha_{0+}\, f)(x) = \begin{cases}
\frac{1}{\Gamma(\alpha)} \int_0^x (x-t)^{\alpha -1}\, f(t)\, dt,\ x>0,& \alpha >0,\\
f(x),\ x>0, & \alpha = 0.
\end{cases}
\end{equation}
The operator $I^\alpha_{0+}$ is well defined, say,  for the functions from the space $L_{loc}(\R_+)$. 
Evidently, it can be interpreted as the Laplace convolution of the functions $f=f(x)$ and $h_\alpha(x) = x^{\alpha -1}/\Gamma(\alpha)$, $x>0$. The convolution theorem for the Laplace transform immediately leads to the well-known formula  for the Laplace transform of the Riemann-Liouville fractional integral $I^\alpha_{0+}$ (\cite{Samko})
\begin{equation}
\label{RLI_L} 
({\mathcal  L}\, I^\alpha_{0+}\, f)(s) = s^{-\alpha}\, ({\mathcal  L}\, f)(s),\ \Re(s) > \max\{s_f,\, 0\}
\end{equation}
that is valid under the condition that the Laplace transform of the function $f$ given by the integral
\begin{equation}
\label{Lap} 
({\mathcal  L}\, f)(s) = \int_0^{+\infty} f(t)\, e^{-st}\, dt
\end{equation}
does exist for $\Re(s) > s_f$. 

Now we rewrite the definition \eqref{nLD} of the $n$th level fractional derivative of order $\alpha,\ 0<\alpha \le 1$ and type $\gamma=(\gamma_1,\dots, \gamma_n)$ for the case of the positive real semi-axis:
\begin{equation}
\label{nLDB_R}
(D^{\alpha,(\gamma)}_{nL_+}\, f)(x)   = \left(\prod_{k=1}^n (I^{\gamma_k}_{0+}\, \frac{d}{dx})\right)\, (I^{n-\alpha-s_n}_{0+}\, f)(x),\ x>0.
\end{equation}

For the appropriate spaces of functions that take into consideration the behavior of the functions at 
$+\infty$ (see \cite{Samko} for more details regarding the case of the Riemann-Liouville fractional integral and derivative on the positive real semi-axis), both the Fundamental Theorem of FC and the projector formula that were  derived for the case of a finite interval remain to be true for the operators $I^{\alpha}_{0+}$ and $D^{\alpha,(\gamma)}_{nL_+}$ defined on the positive real semi-axis. 
The reason is that these formulas are valid point-wise. Thus, for an arbitrary point $x\in \R_+$ we can choose an interval $[0,X]$ with the property $x\in [0,\, X]$. Because the Fundamental Theorem of FC and the projector formula hold true on the whole interval $[0,\, X]$, they are in particular valid at the point $x$.

As already mentioned, for some values of the parameters, the $n$th level fractional derivatives are reduced to the fractional derivatives of the level less than $n$. In the following discussions we restrict ourselves to the truly $n$th level fractional derivatives and suppose that the conditions \eqref{cond_add_n} are satisfied. Then the 2nd Fundamental Theorem of FC (Remark \ref{r_2ndFT}) holds true and we have the representation
$$
(I^{\alpha}_{0+} D^{\alpha,(\gamma)}_{nL_+}\, f)(x) = f(x) - \sum_{k=1}^n \, p_k\, x^{\alpha +s_k -k},
$$
where the coefficients $p_k,\ k=1,\dots,n$ are given by \eqref{p}. 

Under the conditions \eqref{cond_add_n}, the exponents $\sigma_k = \alpha +s_k -k,\ k=1,\dots,n$ satisfy the inequalities $-1 <\sigma_k\le 0$. Thus we can  
apply the Laplace transform to the last formula and get the following equation using the formula \eqref{RLI_L}: 
$$
s^{-\alpha}\, ({\mathcal  L}\,  D^{\alpha,(\gamma)}_{nL_+}\, f)(s) = ({\mathcal  L}\, f)(s) - \sum_{k=1}^n \, p_k\,  \frac{s^{k-\alpha-s_k-1}}{\Gamma(\alpha+s_k-k+1)}.
$$
Now we solve this equation for the Laplace transform of the $n$th level fractional derivative:
$$
({\mathcal  L}\,  D^{\alpha,(\gamma)}_{nL_+}\, f)(s) = s^{\alpha}\, ({\mathcal  L}\, f)(s) -  \sum_{k=1}^n \, p_k\,  \frac{s^{k-s_k-1}}{\Gamma(\alpha+s_k-k+1)}.
$$
Taking into consideration the formula \eqref{p} for the coefficients $p_k,\ k=1,\dots,n$, we arrive at the final formula
\begin{equation}
\label{Lap_D_n}
({\mathcal  L}\,  D^{\alpha,(\gamma)}_{nL_+}\, f)(s) = s^{\alpha}\, ({\mathcal  L}\, f)(s) - \sum_{k=1}^n \, a_k\,   s^{k-s_k-1},
\end{equation}
where the coefficients $a_k,\ k=1,\dots,n$ are determined by the function $f$:
\begin{equation}
\label{a}
a_k= \left( \prod_{i=k+1}^n \left(I^{\gamma_i}_{0+}\, \frac{d}{dx}\right) \, I^{n-\alpha-s_n}_{0+}\, f\right)(0).
\end{equation} 

If some of the conditions \eqref{cond_add_n} does not hold true, the kernel of the $n$th level fractional derivative has a dimension less than $n$ and the corresponding coefficients $a_k$ in the Laplace transform formula \eqref{Lap_D_n} are equal to zero.  

As an example, let us consider the Laplace transform formula for the 2nd level fractional derivative that was derived in \cite{Luc20}:
\begin{equation}
\label{Lap_D}
({\mathcal  L}\,  D^{\alpha,(\gamma_1,\gamma_2)}_{2L_+}\, f)(s) = s^{\alpha}\, ({\mathcal  L}\, f)(s) - a_1\, s^{ - \gamma_1}  - a_2\,   s^{- \gamma_1-\gamma_2+1}
\end{equation}
with
\begin{equation}
\label{a1}
a_1 = \left( I^{\gamma_2}_{0+}\, \frac{d}{dx} \, I^{2-\alpha-\gamma_1-\gamma_2}_{0+}\, f\right)(0), \ 
a_2= \left(  I^{2-\alpha-\gamma_1-\gamma_2}_{0+}\, f\right)(0).
\end{equation} 
For the Hilfer fractional derivative $D^{\alpha,\gamma_1}_{H_+}$ defined on the real positive semi-axis ($\gamma_2 = 1$  or $\alpha +\gamma_1+\gamma_2\le 1$ in  \eqref{nLDB_R}), one of the coefficients $a_1$ (if $\gamma_2 = 1$)  or $a_2$ (if $\alpha +\gamma_1+\gamma_2\le 1$) is zero and we arrive at the known formula
\begin{equation}
\label{Lap_D_H}
({\mathcal  L}\,  D^{\alpha,\gamma_1}_{H_+}\, f)(s) = s^{\alpha}\, ({\mathcal  L}\, f)(s) - \left(  I^{1-\alpha-\gamma_1}_{0+}\, f\right)(0)\,   s^{- \gamma_1}.
\end{equation}
The case $\gamma_1=0$ in \eqref{Lap_D_H} corresponds to the Laplace transform  formula for the Riemann-Liouville fractional derivative
\begin{equation}
\label{Lap_D_RL}
({\mathcal  L}\,  D^{\alpha}_{0+}\, f)(s) = s^{\alpha}\, ({\mathcal  L}\, f)(s) - \left(  I^{1-\alpha}_{0+}\, f\right)(0),
\end{equation}
whereas the case $\gamma_1 = 1- \alpha$ leads to the Laplace transform of the Caputo fractional derivative:
\begin{equation}
\label{Lap_D_C}
({\mathcal  L}\,  D^{\alpha}_{C_+}\, f)(s) = s^{\alpha}\, ({\mathcal  L}\, f)(s) - f(0)\,   s^{\alpha -1}.
\end{equation}

\section{Fractional relaxation equation with the $n$th level fractional derivative}

We start this section with a short discussion of the simplest fractional differential equation with the $n$th level fractional derivative of  order $\alpha,\ 0<\alpha \le 1$, namely, the one in the form
\begin{equation}
\label{de_1}
(D^{\alpha,(\gamma)}_{nL_+}\, y)(x) = 0,\ x>0.
\end{equation}
In fact, we already considered this equation in the previous section because its solution is the kernel of $D^{\alpha,(\gamma)}_{nL_+}$. Depending on the parameters  $\gamma_k, k=1,\dots,n$, the kernel dimension is ranging from $n$ to 1 because the $n$th level fractional derivative can degenerate to the fractional derivatives with the level $n-1$ to 1. In the case, the conditions \eqref{cond_add_n} are satisfied, the kernel of $D^{\alpha,(\gamma)}_{nL_+}$ is $n$-dimensional and  the general solution to the equation \eqref{de_1} is as follows:
\begin{equation}
\label{eq_1}
y(x) = \sum_{k=1}^n \, c_k\, x^{\sigma_k},\ \sigma_k = \alpha +s_k -k,\ x>0,
\end{equation}
$c_k \in \R,\ k=1,\dots,n$ being arbitrary constants. To guarantee uniqueness of solution to the equation \eqref{de_1}, $n$ initial conditions 
\begin{equation}
\label{y-k}
\left( \prod_{i=k+1}^n \left(I^{\gamma_i}_{0+}\, \frac{d}{dx}\right) \, I^{n-\alpha-s_n}_{0+}\, y\right)(0)\ = \ y_k,\ k=1,\dots,n
\end{equation}
are required.
The initial-value problem for the equation \eqref{de_1} with the initial conditions \eqref{y-k} has then the unique solution 
\begin{equation}
\label{sol_1}
y(x) = \sum_{k=1}^n \, \frac{y_k}{\Gamma(\sigma_k+1)}\, x^{\sigma_k},\ x>0.
\end{equation}
This situation is rather unusual for the fractional differential equations with a fractional derivative of the order $\alpha \in ]0,1]$ and has consequences for their possible applications as we will see in the further discussions. 

Now we consider the fractional relaxation equation with the $n$th level fractional derivative and prove the following result:

\begin{The}
\label{t_cm}
The fractional relaxation equation
\begin{equation}
\label{de_2}
(D^{\alpha,(\gamma)}_{nL_+}\, y)(x) = -\lambda\, y(x),\ \lambda >0,\ x>0
\end{equation}
with the initial conditions given by \eqref{y-k} has a unique solution given by the formula
\begin{equation}
\label{sol_4}
y(x) = \sum_{k=1}^n \, y_k\,   x^{\alpha+s_k-k}\, E_{\alpha,\alpha+s_k-k+1}(-\lambda x^\alpha).
\end{equation}
In the case, the initial conditions are non-negative ($y_k\ge 0,\ k=1,\dots,n$ in \eqref{y-k}) and the conditions
\begin{equation}
\label{cond_cm}
k-1 \le s_k,\ k=1,\dots,n.
\end{equation}
hold true, the solution \eqref{sol_4} is completely monotone. 
\end{The}

In the formulation of the theorem, $E_{\alpha,\beta}$ stands for the two-parameters Mittag-Leffler function that is defined by the following convergent series:
\begin{equation}
\label{ML}
E_{\alpha,\beta}(z) = \sum_{k=0}^\infty \frac{z^k}{\Gamma(\alpha\, k + \beta)},\ \alpha >0,\ \beta,z\in \C.
\end{equation}

To prove the theorem, we apply the Laplace transform method. First we do this formally, but then verify that the Laplace transform of the obtained solution does exist. Using the  formula  \eqref{Lap_D_n}, the initial-value problem \eqref{y-k}, \eqref{de_2} can be transformed to the Laplace domain:
\begin{equation}
\label{sol_2}
s^{\alpha}\, ({\mathcal  L}\, y)(s) - \sum_{k=1}^n \, y_k\,   s^{k-s_k-1} \, = \, -\lambda\, ({\mathcal  L}\, y)(s).
\end{equation}
The solution  in the Laplace domain is as follows:
\begin{equation}
\label{sol_3}
({\mathcal  L}\, y)(s) = \sum_{k=1}^n \, y_k\,   \frac{s^{k-s_k-1}}{s^\alpha +\lambda}.
\end{equation}
Now we employ the well-known formula
\begin{equation}
\label{Lap_ML}
({\mathcal  L}\, x^{\beta -1}\, E_{\alpha,\beta} (-\lambda\, x^\alpha)(s) \, = \, \frac{s^{\alpha -\beta}}{s^\alpha +\lambda},\ \beta >0,
\end{equation}
and immediately arrive at the formula \eqref{sol_4} for the solution to the fractional relaxation equation \eqref{de_2} with the initial conditions  \eqref{y-k}.

As already mentioned in the previous section, under the conditions \eqref{cond_add_n}, the exponents $\sigma_k=\alpha+s_k-k,\ k=1,\dots,n$ fulfill the inequalities $-1 < \sigma_k \le 0$. For $\alpha \in ]0,\, 1]$, the Mittag-Leffler function has the following asymptotics as $x\to -\infty$:
\begin{equation}
\label{ML_a}
E_{\alpha,\beta}(x) = -\frac{x^{-1}}{\Gamma(\beta-\alpha)} + O(x^{-2}),\ x\to -\infty.
\end{equation}
These both facts ensure that the Laplace transform of the function at the right-hand side of \eqref{sol_4} does exists and thus it is indeed the unique solution to the initial-value problem \eqref{y-k}, \eqref{de_2}. 

Now let us prove that the solution \eqref{sol_4} is a completely monotone function provided the conditions \eqref{cond_cm} hold true and the initial conditions are all non-negative. 

For the reader's convenience, we recall that a  non-negative function $ \phi:\R_+\to \R$  is called  completely
monotone if it is from $C^{\infty}(\R_+)$  and
$(-1)^n \phi^{(n)}(x)\geq 0$ for all $n \in \N$ and
$x\in \R_+$.

The Mittag-Leffler function $f(x)=E_{\alpha,\beta}(-x)$ is completely monotone if and only if  $0<\alpha \le 1$ and $\alpha \le \beta$ (\cite{Sch}). The power law function $g(x) = \lambda\, x^{\alpha},\ 0< \alpha \le 1,\ 0<\lambda$ is a Bernstein function because its derivative $g^\prime (x) = \lambda\alpha\, x^{\alpha-1}$ is completely monotone. A composition of a completely monotone function and a Bernstein function is completely monotone (\cite{SSV}). Then the  function $E_{\alpha,\beta}(-\lambda x^{\alpha})$ is completely monotone for $0<\alpha \le 1$, $\alpha \le \beta$, and $0<\lambda$. The product of two completely monotone functions is again a completely monotone function (\cite{SSV}) that leads to the complete monotonicity of the function 
\begin{equation}
\label{cm}
h_{\alpha,\beta,\gamma,\lambda}(x) := x^{\gamma-1}\, E_{\alpha,\beta}(-\lambda x^{\alpha})
\end{equation}
under the conditions
\begin{equation}
\label{cm_c}
0<\alpha \le 1,\ \alpha \le \beta,\ 0 < \gamma \le 1,\ 0<\lambda.
\end{equation}
The functions $y_k(x) = x^{\alpha+s_k-k}\, E_{\alpha,\alpha+s_k-k+1}(-\lambda x^\alpha),\ k=1,\dots,n$ from the solution formula \eqref{sol_4} have the form \eqref{cm}. Moreover, their parameters $\alpha = \alpha$, $\beta = \alpha+s_k-k+1$, and $\gamma = \alpha+s_k-k+1$ fulfill the conditions \eqref{cm_c} because of the
conditions \eqref{gamma}, \eqref{cond_add_n}, and \eqref{cond_cm} we posed on the order $\alpha$ and the type $\gamma =(\gamma_1,\dots,\gamma_n)$ of the $n$th level fractional derivative. Thus, the functions $y_k(x)$ are all completely monotone as well as their linear combination with non-negative coefficients that builds the solution \eqref{sol_4}. The proof of the theorem is completed.

\begin{Rem}
Taking into account the formula \eqref{ML_a} for the asymptotics of the Mittag-Leffler function, the behavior of the solution \eqref{sol_4} as $x\to +\infty$ has the following form:
\begin{equation}
\label{asymp}
y(x) \sim  \sum_{k=1}^n \, d_k\, x^{\alpha_k},\ \alpha_k = s_k-k,\ x\to +\infty,
\end{equation}
where the coefficients $d_k$ depend both on the initial values $y_k$, the order $\alpha$, and the type $\gamma = (\gamma_1,\dots,\gamma_n)$. Thus, the fractional relaxation equation \eqref{de_2} can be employed to model the relaxation processes with the asymptotic behavior of type \eqref{asymp}. Moreover, the free parameters $\gamma_1,\dots,\gamma_n$ can be used for optimal fitting of the measurements data for a concrete relaxation process with a power law asymptotics. 
\end{Rem} 

\begin{Rem}
In \cite{DN}, uniqueness and existence of solutions to the Cauchy problems for the linear and non-linear fractional differential equations with the Djrbashian-Nersesian operators similar to the $n$th level fractional derivatives $D^{\alpha,(\gamma)}_{nL}$ (see Remark \ref{r_dn}) were addressed. In particular, an analogy of the formula \eqref{sol_4} on a finite interval was deduced  using the method of power series. However, no analysis of the solution properties including their complete monotonicity was presented there. 

It is worth mentioning that in \cite{D_book}, the eigenfunctions and the associated functions for some boundary value problems for the special equations containing the Djrbashian-Nersesian operators were constructed in explicit form. Moreover, these functions were interpreted as the biorthogonal systems of vector functions and used for the interpolation expansions for some Hilbert spaces of entire functions.   
\end{Rem}

In the rest of this section, we illustrate Theorem \ref{t_cm} on the case of the fractional relaxation equation with the truly 2nd level derivative (the conditions \eqref{cond_add_n} with $n=2$ hold true)
\begin{equation}
\label{de_3}
(D^{\alpha,(\gamma_1,\gamma_2)}_{2L_+}\, y)(x) = -\lambda\, y(x),\ \lambda >0,\ x>0
\end{equation}
and with the initial conditions 
\begin{equation}
\label{y12}
\left( I^{\gamma_2}_{0+}\, \frac{d}{dx} \, I^{2-\alpha-\gamma_1-\gamma_2}_{0+}\, f\right)(0)\, = \, y_1, \ 
\left(  I^{2-\alpha-\gamma_1-\gamma_2}_{0+}\, f\right)(0) \, = \, y_2.
\end{equation} 
According to Theorem \ref{t_cm}, its solution 
\begin{equation}
\label{sol_5}
y(x) = y_1\,   x^{\alpha+\gamma_1-1}\, E_{\alpha,\alpha+\gamma_1}(-\lambda x^\alpha) +
y_2\,   x^{\alpha+\gamma_1+\gamma_2-2}\, E_{\alpha,\alpha+\gamma_1+\gamma_2-1}(-\lambda x^\alpha)
\end{equation} 
is completely monotone if the initial conditions $y_1,\, y_2$ are non-negative and the following restrictions on the order $\alpha$ and type $(\gamma_1,\, \gamma_2)$ of the 2nd level fractional derivative hold true:
\begin{equation}
\label{cond_2l}
0\le \gamma_1\le 1-\alpha,\ 0\le \gamma_2 < 1,\ 1\le \gamma_1 + \gamma_2 \le 2- \alpha.
\end{equation}
The points of the $(\gamma_1,\, \gamma_2)$-plane that satisfy the conditions \eqref{cond_2l} are graphically represented in the plot of Figure 1.  
\begin{figure}[htbp]
\label{fig}
	\centering
		\includegraphics[scale=0.8]{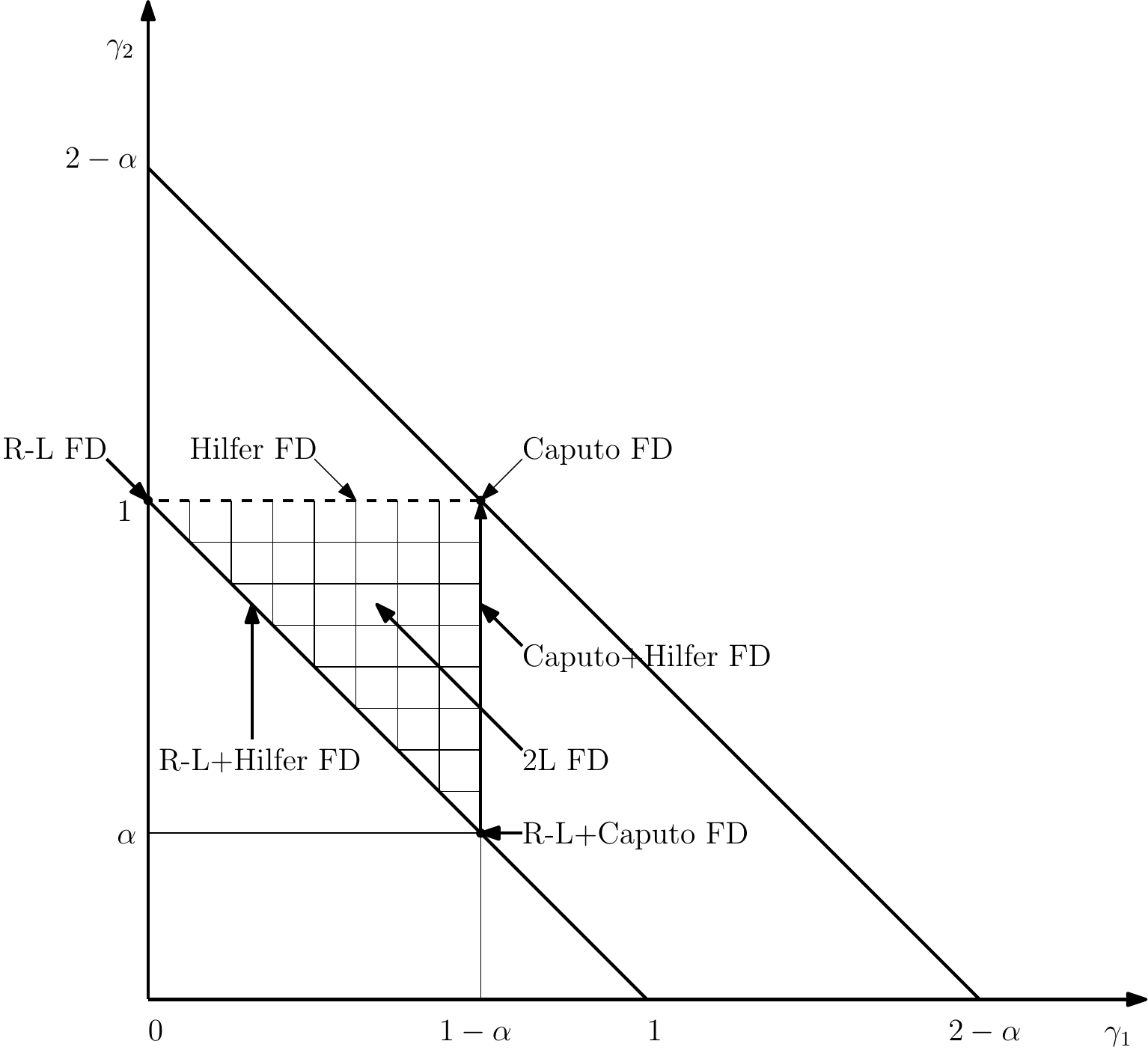}
	\caption{2nd level fractional derivatives of order $\alpha\in ]0,1]$ and type $(\gamma_1,\, \gamma_2)$}
	\label{fig:CT-bild}
\end{figure} 
They build a triangle and in what follows we shortly discuss the particular cases of the fractional relaxation equation \eqref{de_2} that correspond to the vertexes and edges of this triangle. As already mentioned in the previous section, the 2nd level fractional derivatives with $\gamma_2 = 1$ (upper edge of the triangle) are reduced to the 1st level derivatives (the Riemann-Liouville, the Caputo, and the Hilfer fractional derivatives). Still, it is instructive to include them into considerations. 

We start with the vertexes and write down the types of the corresponding fractional derivatives, their form, and the solutions to the fractional relaxation equation \eqref{de_2} with these derivatives:

\begin{itemize}

\item The Riemann-Liouville fractional derivative: 

\begin{itemize}

\item $\gamma_1 =0,\ \gamma_2 =1$, 
\item[] $ $
\item $D^{\alpha,(\gamma_1,\gamma_2)}_{2L_+} = \frac{d}{dx}\, I^{1-\alpha}_{0+}$, 
\item[]$ $
\item $y(x) = y_1\, x^{\alpha -1}\, E_{\alpha,\alpha}(-\lambda x^{\alpha})$. 

\end{itemize}

\item  The Caputo fractional derivative: 

\begin{itemize}

\item $\gamma_1 =1-\alpha,\ \gamma_2 =1$, 
\item[]$ $
\item $D^{\alpha,(\gamma_1,\gamma_2)}_{2L_+} = I^{1-\alpha}_{0+}\, \frac{d}{dx}$, 
\item[] $ $
\item $y(x) = y_1\,  E_{\alpha,1}(-\lambda x^{\alpha})$. 

\end{itemize}

\item A truly 2nd level fractional derivative:

\begin{itemize}

\item  $\gamma_1 =1-\alpha,\ \gamma_2 =\alpha$, 
\item[] $ $
\item $D^{\alpha,(\gamma_1,\gamma_2)}_{2L_+} = I^{1-\alpha}_{0+}\, \frac{d}{dx}\, I^{\alpha}_{0+}\, \frac{d}{dx}\,I^{1-\alpha}_{0+}$, 
\item[] $ $
\item $y(x) = y_1\,  E_{\alpha,1}(-\lambda x^{\alpha}) + y_2\, x^{\alpha -1}\, E_{\alpha,\alpha}(-\lambda x^{\alpha}) $. 

\end{itemize}

\end{itemize}

Whereas the first two cases (the fractional relaxation equations with the Riemann-Liouville and the Caputo fractional derivatives) are well-known, the third case seems to be new. The solution to the corresponding relaxation equation is a linear combination of the solutions to the fractional relaxation equations with the Riemann-Liouville and with the Caputo fractional derivatives. 

Now let us inspect the edges of the triangle from Figure 1. 

\begin{itemize}

\item The Hilfer fractional derivative: 

\begin{itemize}

\item $0\le \gamma_1 \le 1-\alpha,\ \gamma_2 =1$, 
\item[] $ $
\item $D^{\alpha,(\gamma_1,\gamma_2)}_{2L_+} = I^{\gamma_1}_{0+}\, \frac{d}{dx}\, I^{1-\alpha-\gamma_1}_{0+}$, 
\item[] $ $
\item $y(x) = y_1\, x^{\alpha+\gamma_1 -1}\, E_{\alpha,\alpha+\gamma_1}(-\lambda x^{\alpha})$.

\end{itemize}

\item A truly 2nd level fractional derivative: 

\begin{itemize}

\item $\gamma_1 =1-\alpha,\ \alpha < \gamma_2 <1$, 
\item[]$ $
\item $D^{\alpha,(\gamma_1,\gamma_2)}_{2L_+} = I^{1-\alpha}_{0+}\, \frac{d}{dx}\, I^{\gamma_2}_{0+}\, \frac{d}{dx}\,I^{1-\gamma_2}_{0+}$, 
\item[] $ $
\item $y(x) = y_1\,  E_{\alpha,1}(-\lambda x^{\alpha}) + y_2\, x^{\gamma_2 -1}\, E_{\alpha,\gamma_2}(-\lambda x^{\alpha}) $. 

\end{itemize}

The solution  to the relaxation equation with the 2nd level fractional derivative $D^{\alpha,(1-\alpha,\gamma_2)}_{2L_+}$ is a linear combination of the solutions to the fractional relaxation equations with the Caputo fractional derivative and with the Hilfer fractional derivative with the type $\gamma_1 = \gamma_2 - \alpha$. 

\item A truly 2nd level fractional derivative: 

\begin{itemize}

\item $0<\gamma_1 < 1-\alpha,\ \gamma_2 = 1-\gamma_1$, 
\item[] $ $
\item $D^{\alpha,(\gamma_1,\gamma_2)}_{2L_+} = I^{\gamma_1}_{0+}\, \frac{d}{dx}\, I^{1-\gamma_1}_{0+}\, \frac{d}{dx}\,I^{1-\alpha}_{0+}$, 
\item[] $ $
\item $y(x) = y_1\, x^{\alpha+\gamma_1 -1}\, E_{\alpha,\alpha+\gamma_1}(-\lambda x^{\alpha})  +y_2\,  x^{\alpha -1}\, E_{\alpha,\alpha}(-\lambda x^{\alpha}) $. 

\end{itemize}

The solution  to the relaxation equation with the 2nd level fractional derivative $D^{\alpha,(\gamma_1,1-\gamma_1)}_{2L_+}$ is a linear combination of the solutions to the fractional relaxation equations with the Riemann-Liouville fractional derivative and with the Hilfer fractional derivative with the type $\gamma_1$. 

\end{itemize}

In all other cases the solution to the fractional relaxation equation \eqref{de_3} is given by the formula \eqref{sol_5}. It is worth mentioning that \eqref{sol_5} can be interpreted as a linear combination of solutions to the fractional relaxation equations with the Hilfer fractional derivatives of order $\alpha$ and with the types $\gamma_1$ and $\gamma_1+\gamma_2-1$, respectively. 

Finally let us mention that the 2nd Fundamental Theorem of FC for the $n$th order fractional derivative (Remark \ref{r_2ndFT}) can be used for analysis of more complicated and even nonlinear fractional differential equations. Say, the fractional differential equation
\begin{equation}
\label{eq-nl}
(D^{\alpha,(\gamma)}_{nL_+}\, y)(x) \, =\, F(x,y(x)),\ x>0,\ F:\, \R_+\times \R \to \R
\end{equation}
subject to the initial conditions in form \eqref{y-k} can be transformed to the following Volterra-type integral equation of the second kind by applying the Riemann-Liouville fractional integral to the equation \eqref{eq-nl} and by using the formula \eqref{2ndFT}:
\begin{equation}
\label{eq-nl_sol}
y(x) \, = \, (I^{\alpha}_{0+}\,  F(x,y(x)))(x) \, + \, \sum_{k=1}^n \, \frac{y_k}{\Gamma(\sigma_k +1)}\, x^{\sigma_k}, \ \sigma_k = \alpha +s_k -k.
\end{equation}
The integral equation \eqref{eq-nl_sol} can be analyzed by the standard method of the fix point iterations. This problem will be considered elsewhere.   





\end{document}